\documentclass{article}

\begin{document}

\title{K\"ahler metrics with cone singularities along a divisor}
\author{S. K. Donaldson}
\maketitle

\newtheorem{conj}{Conjecture}
\newtheorem{thm}{Theorem}
\newtheorem{prop}{Proposition}
\newtheorem{cor}{Corollary}
\newtheorem{lem}{Lemma}
\newcommand{\cD}{{\cal D}}
\newcommand{\osigma}{\overline{\sigma}}
\newcommand{\grad}{{\rm grad}}
\newcommand{\bR}{{\bf R}}
\newcommand{\bC}{{\bf C}}
\newcommand{\bP}{{\bf P}}
\newcommand{\db}{\overline{\partial}}
\newcommand{\ut}{\underline{t}}
\newcommand{\up}{\underline{p}}
\newcommand{\Vol}{{\rm Vol}}
\newcommand{\Ric}{{\rm Ric}}
\newcommand{\dbd}{i\partial \overline{\partial}}
\newcommand{\ozeta}{\overline{\zeta}}
\section{Introduction}
Let $D$ be a smooth divisor in a complex manifold $X$. In this paper we study K\"ahler metrics on $X\setminus D$ with  cone singularities of cone angle $2\pi \beta$  transverse to $D$, where $0<\beta<1$. The case we have primarily in mind is  when $X$ is a Fano manifold, $D$ is an anticanonical divisor and the metrics are K\"ahler-Einstein; the motivation being the hope that one can study the existence problem for {\it smooth} K\"ahler-Einstein metrics on $X$
(as a limit when  $\beta$ tends to $1$) by deforming the cone angle. This can be seen as a variant of the standard \lq\lq continuity method''. We will make some more remarks about this programme in Section 6 but it is clear that, at the best, a substantial amount of work will be needed to  carry this through--- adapting much of the standard theory  to the case of cone singularities.
This paper is merely a first step along this  road. Our goal is  to set up a linear theory and apply it to the problem of deforming the cone angle (Theorem 2 below). In further papers with X-X Chen, we will study more advanced and sophisticated questions. 

There are several precedents for this line of work. First and foremost, singular metrics of this kind have been considered before by Jeffres \cite{kn:TJ},\cite{kn:TJ2} and Mazzeo \cite{kn:Mazzeo}. Some applications to algebraic geometry are outlined by Tian in \cite{kn:GT}.  Mazzeo considers the case of negative first Chern class, but this makes no difference in the elementary foundational questions we consider here (until Section 6). As in this paper, Mazzeo's main emphasis is on the linear theory, and he outlines an approach using the \lq\lq edge calculus''. However this assumes some specialised background, some complications with the choice of function spaces are reported and \cite{kn:Mazzeo} does not give quite enough detail for those not expert in the techniques to easily fill in the proofs. Thus  we have decided  to make a fresh start here on the analysis, using elementary methods. This means that we are very probably re-deriving many results that are well-known to experts, and our conclusions are entirely consistent with those described by Mazzeo. It is very likely that the edge calculus, or similar technology,  will be important in developing more refined analytical results. 

A second precedent occurs in the study of 3-dimensional hyperbolic manifolds. Here again one can consider metrics with cone singularities transverse to a knot. A strategy, similar  to ours in K\"ahler geometry, for constructing nonsingular hyperbolic metrics via deformation of the cone angle was proposed by Thurston and there are a number of papers in the literature developing the theory and the relevant analysis (\cite{kn:HK}, \cite{kn:MM}, \cite{kn:MW}\ for example).   A third precedent occurs in gauge theory and the work of Kronheimer, Mrowka and others on connections with codimension-2 singularities\cite{kn:KM}.
In the case when the underlying manifold is complex, this is related to the theory of holomorphic bundles with parabolic structures \cite{kn:B} and  there are some closer parallels with our situation.

The general scheme of this paper mimics the development of standard theory for smooth manifolds. We begin by considering a \lq\lq flat model'' for a cone singularity and in  Section 2 we obtain an estimate in H\"older spaces for the Laplace operator, as in the usual Schauder theory. This depends on certain properties of the Green's function which are derived in Section 3, using Bessel functions and classical methods. In Section 4 we introduce complex structures, considering first a  flat model and then a general class of singular metrics on a pair $(X,D)$. What we achieve is roughly, a  parallel to the standard theory of H\"older continuous K\"ahler metrics. This degree of regularity  suffices to give a Fredholm theory linearising the K\"ahler-Einstein equation, and in particular we can proceed to study the problem of deforming the cone angle.  Naturally we expect that it will be possible to say much more about the local structure of these solutions but we mainly leave this for future papers.  Sections 5 and 6 are intended to provide context.  In Section 5 we use the Gibbons-Hawking construction, combined with our study of the Green's function, to produce certain almost-explicit Ricci-flat metrics with cone singularities, analogous to ALE spaces in the usual theory. In Section 6 we outline what one might expect when $X$ is the complex projective plane  blown up at one or two points---when no smooth K\"ahler-Einstein metrics exist---and  discuss connections with work of Szekelyhidi and Li.

The author is grateful to Xiu-Xiong Chen, Mark Haskins and Jared Wunsch for discussions related to this work.

\section{A Schauder estimate}

For $\alpha\in (0,1)$ and for a function $f$ on $\bR^{m}$ we define
$$    [f]_{\alpha}= \sup_{p,q} \frac{\vert f(p)-f(q)}{\vert p-q\vert^{\alpha}}, $$
where $\sup=\infty$ is allowed. Write $\bR^{m}=\bR^{2}\times \bR^{m-2}$ and   let $S= \{0\}\times \bR^{m-2}$. Take polar co-ordinates
$r,\theta$ on $\bR^{2}$ and standard co-ordinates $s_{i}$ on $\bR^{m-2}$.
Fix $\beta \in (0,1)$ and consider the singular metric 
\begin{equation} g= dr^{2} + \beta^{2} r^{2} d\theta^{2} + \sum ds_{i}^{2}. \end{equation}
This is the {\it standard cone metric} with cone angle $2\pi \beta$ and a singularity along  $S$. 
We want to consider the Green's operator of the Laplacian $\Delta=\Delta_{g}$. To fix a definition, let $H$ be the be the completion of $C^{\infty}_{c}$ under the Dirichlet norm $\Vert \nabla f \Vert_{L^{2}}$. Since the metric $g$ is uniformly equivalent to the standard Euclidean  one, we get the same space $H$ using either metric. The Sobolev inequality implies that for $q=2m/(m+2)$ and any  $\rho\in L^{q} $ the linear form  
$$  f\mapsto \int f \rho, $$
is bounded with respect to the $H$ norm, so there is a unique $G\rho\in H$ such that
$$  \int f \rho = \int (\nabla f, \nabla G \rho)_{g}, $$ which is to say that $\phi=G\rho$ is a weak solution of the equation $\Delta_{g}\phi=\rho$.
Thus we define a linear map $G:L^{q}\rightarrow H$.

\begin{prop}
There is a locally-integrable kernel function $G(x,y)$ such that
$$   G\rho(x)= \int G(x,y) \rho(y) dy, $$
for $\rho\in C^{\infty}_{c}$. The function $G(x,y)$ is smooth away from the diagonal and points $x,y\in S$.
\end{prop}

This follows from standard theory, but in the next section we will give an \lq\lq explicit'' formula for $G$.

Let $D$ be one of the differential operators
$$  \frac{\partial^{2}}{\partial s_{i}\partial s_{j}} \ \,\ \ \frac{\partial^{2}}{\partial r \partial s_{i}} \ \ , \ \ \frac{1}{r} \frac{\partial^{2}}{\partial \theta \partial s_{i}}. $$

We define $T=D \circ G$. Let $\mu=\beta^{-1}-1$. The main result of this section is

\begin{thm}
Fix $\alpha$ with $0<\alpha< \mu$. Then there is a constant $C$ depending on $\beta, n,\alpha$ such that for all functions $\rho\in C^{\infty}_{c}(\bR^{m})$ we have
$$  [T\rho]_{\alpha} \leq C [\rho]_{\alpha}. $$
\end{thm}

\

(The statement should be interpreted as including the assertion that $T\rho$ is continuous, so its value at each point is defined.)

{\bf Note} When we refer to the distance $d(x,y)=\vert x-y\vert$ between points in $\bR^{m}$ we always mean the standard Euclidean distance. However this is uniformly equivalent to the distance defined by the singular metric.

\

The proof of the Theorem uses an integral representation for $T$. Let $K(x,y)= D_{x}G(x,y)$ where the notation means that the differentiation is applied to the first variable. Then we have
\begin{prop}
If $\rho\in C^{\infty}_{c}$ and $\rho(x_{0})=0$ for some $x_{0}\in \bR^{m}$ then $K(x_{0}, \ )\rho(\ )$ is integrable and
$$  (T\rho)(x_{0})= \int K(x_{0}, y) \rho(y) dy. $$
\end{prop}

Of course the subtlety is that if $\rho$ does not vanish at $x_{0}$ then $K(x_{0}, \ ) \rho(\ )$ is not integrable and the formula has to be interpreted as a singular integral, but we will not need to use this approach. What we do need is some more detailed information about the kernel $K$, summarised in the next Proposition. We write $\pi:\bR^{2}\times \bR^{m-2}\rightarrow \bR^{2}$ for the projection map. 

\begin{prop}
There are $\kappa_{1}, \kappa_{2}, \kappa_{3},\kappa_{4}$ with the following properties. 

\begin{itemize}\item If $\vert z\vert=1$ then  $$\vert K(0,z) \vert \leq \kappa_{1}$$
\item If $\vert z\vert=1$ then $$ \vert K(w_{1}, z)-K(w_{2}, z) \vert \leq \kappa_{2} \vert w_{1}-w_{2}\vert^{\mu}$$
for any $w_{1}, w_{2}$ with $\vert w_{i}\vert \leq 1/2$
\item If $\vert z\vert=1$ and  $\vert \pi(z)\vert \geq 1/2$ then
$$  \vert K(z,w)\vert\leq \kappa_{3} \vert z-w\vert^{-n}, $$
for $w$ with $\vert w\vert \leq 5$.
\item If $\vert z\vert=1$ and $\vert \pi(z)\vert \geq 1/2$ then 
$$ \vert (\nabla K)(z,w) \vert \leq \kappa_{4} \vert z-w\vert^{-n-1}, $$
for $w$ with $\vert w\vert \leq 5$. Here the derivative $\nabla K$ is taken with respect to the first variable.
\end{itemize}
\end{prop}

Propositions 2 and 3 will be established in the next section but now, assuming them, we go on to the proof of Theorem 1. This is a variant of the standard proof of the Schauder estimate for the ordinary Laplace operator.  

\

What is crucially important is  that $T$ commutes with  dilations. Thus, given $\lambda>0$ and a function $\rho$ on $\bR^{n}$ we define $\rho_{\lambda}(x)= \rho(\lambda^{-1} x)$ and we have $(T\rho)_{\lambda}= T(\rho_{\lambda})$.   This implies that
\begin{equation}  K(\lambda x,\lambda y)= \lambda^{-n} K(x,y).  \end{equation}

Note also that the H\"older seminorm scales by dilation as $[f_{\lambda}]_{\alpha}= \lambda^{-\alpha} [f]_{\alpha}$ so our problem is scale invariant.

Fix a smooth function $\psi$ supported  in the unit ball, with $\Delta_{g} \psi$ and $D \psi$ both smooth and with  $\Delta_{g}\psi=1$  on the $\delta$-ball for some fixed $\delta>0$.  For example we can take $\psi=a(r) b(s)$ where $a(r)$ is equal to $1$ for small $r$ and $b$ is a suitable function of $s$. Set $\chi=\Delta_{g}\psi$ so $\chi$ has compact support, is equal to $1$ on the $\delta$-ball and $T\chi=D \psi$ is smooth.  We write
$[\chi]_{\alpha}=c_{0}$, $[T\chi]_{\alpha}=c_{1}$.

\

By scale invariance and linearity, it suffices to show that if $\rho\in C^{\infty}_{c}$ has $[\rho]_{\alpha}=1$ and if $x_{1}, x_{2}\in \bR^{m}$ with $\vert x_{1}-x_{2}\vert =1$ then
$  \vert \rho(x_{1})-\rho(x_{2})\vert \leq C$.  Let $d=\min(\vert \pi(x_{1})\vert, \vert \pi(x_{2})\vert)$.  We consider two cases:
Case A, when $d\leq 2$, and 
 Case B, when $d>2$.

\

{\bf Case A.} Let $x'_{1}, x'_{2}$ be the projections of $x_{1}, x_{2}$ to $S$. Then we can write
$$  T\rho(x_{1})-T\rho(x_{2})=\left( T\rho(x_{1})- T\rho(x'_{1})\right) +\left( T\rho(x'_{1})- T\rho(x'_{2})\right)+\left( T\rho(x'_{2})- T\rho(x_{2})\right)$$
and $\vert x_{1}- x'_{1}\vert, \vert x'_{1}-x'_{2}\vert, \vert x'_{2}- x_{2}\vert$ are all bounded by $3$. Using this, and translation and scale invariance, it suffices to consider two sub-cases

{\bf Sub-case A1}\ \ \  $x_{1}= 0, \vert x_{2}\vert =1,\  x_{2}\in S$.

{\bf Sub-case A2}\ \ \  $x_{1}= 0, \vert x_{2}\vert =1, \ x'_{2}=0$. (That is, $x_{2}$ lies in $\bR^{2}\times \{0\}$. )

 But to begin with the same discussion applies to either sub-case. We define $\sigma_{0}= \rho(x_{2}) \chi_{\lambda}$ where
$$  \lambda= \max(\delta^{-1}, \vert \rho(x_{2})^{1/\alpha}\vert).$$ 
We also define 
$$  \sigma_{1}= (\rho(0)-\rho(x_{2})) \chi. $$
Then $[\sigma_{0}]_{\alpha}, [\sigma_{1}]_{\alpha}\leq c_{0}$ and 
$[T\sigma_{0}]_{\alpha}, [T\sigma_{1}]_{\alpha}\leq c_{1}$, using the fact that
$\vert \rho(x_{2})-\rho(0)\vert \leq [\rho]_{\alpha}=1$, by hypothesis. Now set $\rho'=\rho-\sigma_{0}-\sigma_{1}$. Thus  $\rho'$ vanishes at $0$ and $x_{2}$ and we have
$$   [\rho']_{\alpha} \leq [\rho]_{\alpha} +2c_{0}=1+2c_{0}\ \ , \ \ \vert T\rho(x_{2})- T\rho(0)\vert \leq \vert T\rho'(x_{2})- T\rho'(0) \vert + 2 c_{1}. $$
This means that, simplifying notation, we can reduce to the situation where
$\rho$ vanishes at $x_{2}$ and $0$. 
Thus, in this situation,  we want to estimate
$$   \int K(0,y) \rho(y) dy - \int K(x_{2},y) \rho(y) dy $$
which is dominated by

$$ I= \int \vert K(x_{2}, y)-K(0,y) \vert \ \vert \rho(y) \vert dy. $$
Consider the contribution from the region $\vert y\vert \geq 2$. By the homogeneity we have
$$  K(x_{2}, y)-K(0,y) = \vert y\vert^{-n}\left(  K\left(\frac{x_{2}}{\vert y\vert}, \frac{y}{\vert y\vert }\right)-K\left(0, \frac{y}{\vert y\vert}\right) \right), $$
and the second item in Proposition 3 gives
$$  \vert K(x_{2}, y)-K(0,y) \vert \leq \kappa_{2} \vert y\vert^{-\mu-n}. $$
Then we get a bound on the contribution to $I$ from $\{\vert y\vert \geq 2\}$ in the form
$$  \int_{2}^{\infty} \kappa_{2} R^{-\mu-n}  R^{\alpha} R^{n-1} dR, $$
which is finite since  $\alpha<\mu$.

Next we have to  estimate the contribution to $I$ from $\{\vert y\vert \leq 2\}$. 
First we consider
$$ I_{1}= \int_{\vert y\vert \leq 2}\vert K(0,y)  \rho(y)\vert dy. $$
By the homogeneity and the first item of Proposition 3 we have
$$  \vert K(0,y) \vert \leq \kappa_{1} \vert y\vert^{-n}, $$
and $\vert \rho(y)\vert \leq \vert y \vert^{\alpha}$ so
$$  I_{1}\leq \kappa_{1} \int_{\vert y\vert \leq 2} \vert y \vert^{-n+\alpha}$$
which is finite.
The final step is to estimate
$$  I_{2}= \int_{\vert y\vert \leq 2} \vert K(x_{2}, y)\rho(y) \vert dy. $$
This is where we use different arguments in the two sub-cases. In sub-case A1, when $x_{2}$ lies in $S$, the estimate is just the same as for $I_{1}$ above, using translation invariance in the $\bR^{m-2}$ factor.  In sub-case A2, when $x_{2}$ is in the orthogonal complement of $S$, we use the third item of Proposition 3 to get
$$ \vert K(x_{2},y)  \vert \leq \kappa_{3} \vert y-x_{2}\vert^{-n}$$
when $\vert y\vert \leq 2$ and so
$$  \vert K(x_{2},y)  \vert \vert\rho(y) \vert \leq \kappa_{3} \vert y-x_{2}\vert^{\alpha-n}$$
and we can proceed as before.  This completes the proof for  Case A. 

\

\

{\bf Case B}

\

Recall  that  we  have $x_{1}, x_{2}$ with $\vert x_{1}-x_{2} \vert =1$ and $\vert \pi(x_{i})\vert >2$. Set $\lambda=\max(\vert x_{1}\vert, \vert x_{2} \vert, \vert \rho(x_{2})\vert ^{1/\alpha})$ and define $\sigma_{0}= \rho(x_{2})\chi_{\lambda}$. Then $\sigma_{0}(x_{1})=\sigma_{0}(x_{2})=\rho(x_{2})$ and we have bounds on $[\sigma_{0}]_{\alpha}, [T\sigma_{0}]_{\alpha}$ as before. It is clear that we can choose a function $\tilde{\psi}$, supported in the unit ball centred at $x_{1}$, with $\Delta_{g} \tilde{\psi}$ equal to $1$ in a small neighbourhood of $x_{1}$ and in such a way that $[\tilde{\psi}]_{\alpha}, [\Delta \tilde{\psi}]_{\alpha}$ are bounded by fixed constants, independent of $x_{1}$ provided only that  $\vert \pi (x_{1})\vert>2$ (that is, $x_{1}$ stays well away from the singular set). Then we put $\tilde{\chi}=\Delta_{g} \tilde{\psi}$ and $$  \rho'= \rho- ( \sigma_{0}+ (\rho(x_{1}-\rho(x_{2})) \tilde{\chi}. $$
Arguing just as before, we are reduced to the situation where $\rho(x_{1})=\rho(x_{2})=0$.

Now we can obviously suppose that $x_{1}$ is the point closest to $S$ and by translation we can suppose that $\vert x_{1}\vert =d$.   We have to estimate the integral $I$, as before. We consider the contribution from three regions
\begin{itemize}
\item Points $y$ with $\vert y\vert>2d$. This goes just as in Case A, using the second item in Proposition 3, and rescaling.
\item Points $y$ with $\vert y-x_{1}\vert \leq 2$. This goes just as before using the third item in Proposition 3 and rescaling.
\item Points $y$ with $\vert y\vert \leq 2d $ and $\vert y-x_{1}\vert>2$.
\end{itemize}
Here we use the fourth item in Proposition 3. Set $z_{i}= x_{i}/d$ and $w=y/d$.
The fourth item in Proposition 3 gives a bound on the derivative of $K(z,w)$ with respect to $z$ for all points $z$ on the segment joining $z_{1}, z_{2}$. For such points the distance $\vert z-w\vert $ is comparable to $\vert z_{1}-w\vert $ so, integrating the bound gives
$$  \vert K(z_{1}, w)-K(z_{2}, w) \vert \leq \kappa d^{-1} \vert w-z_{1} \vert^{-n-1}. $$
since the distance $\vert z_{1}-z_{2} \vert$ is $d^{-1}$. Scaling back we have
$$  \vert K(x_{1}, y)- K(x_{2}, y) \vert \leq \kappa \vert y-x_{1} \vert^{-n-1}. $$
Now we can bound the contribution to $I$ from this region by
$$  \kappa \int_{2}^{3d} R^{-n-1} R^{\alpha} R^{n-1} dR <\infty, $$
where $R=\vert y-x_{1}\vert$.

\section{Representation of the Green's functions by Bessel functions}

Write $c=\beta^{-1}$ and consider the map
$\iota:\bR^{2}\times \bR^{m-2}\rightarrow \bR^{2}\times \bR\times \bR^{m-2}$ defined by
$\iota(r\cos \theta, r\sin\theta,s)= (r^{c}\cos\theta,r^{c}\sin \theta, r^{2}, s)$. For an open subset $\Omega\subset \bR^{m}$ we say that function $f$ on $\Omega$ is $\beta$-smooth if each point of $\Omega$ has a neighbourhood $N\subset \Omega$ such that the restriction of $f$ to $N$ is the composite of $\iota$ and a smooth function in the ordinary sense on a neighbourhood of $\iota(N)$.  We  define the notion of convergence of $\beta$-smooth functions similarly. For fixed $y$ write $\Gamma_{y}=G(\ ,y)$. Then we have
\begin{prop}
If $y$ is not in $\Omega$ then $\Gamma_{y}$ is $\beta$-smooth on $\Omega$ and $\Gamma_{y}$ varies continuously with $y$, with respect to the topology of $\beta$-smooth functions on $\Omega$. 
\end{prop}

If the point $y$ is not in $S$ then   we can identify the metric $g$ in a neighbourhood of $y$ with the usual Euclidean metric and it follows from standard theory that $\Gamma_{y}$ differs from the usual Newton potential by a smooth (in fact harmonic) function. It is straightforward to deduce from   Proposition 4, this observation, and the symmetries and scaling behaviour of the Green's function that our kernel $K$ satisfies the criteria stated in Proposition 3.  Likewise for the proof of Proposition 2.  The main point of interest is the second item of Proposition 3: this is the only place where the number $\mu$, and hence the restriction on the range of the H\"older exponent, appears. The derivative of the map
$$   (r\cos \theta, r\sin \theta)\mapsto (r^{c} \cos \theta, r^{c} \sin \theta) $$ is H\"older continuous with exponent $\mu$. Then the chain rule shows that   for a $\beta$-smooth function $f$ the derivatives $\frac{\partial f}{\partial r}$ and $r^{-1} \frac{\partial f}{\partial \theta}$ are H\"older continuous with this exponent. It follows that, for each choice of differential operator $D$, the derivative $D\Gamma_{y}$ is $C^{,\mu}$ near the singular set (the derivatives in the $s_{i}$ variables being harmless).

 Granted the assertions above, we will  focus for the rest of this Section on the proof of Proposition 4. We achieve this by showing that the Green's function has a \lq\lq polyhomogeneous expansion'' around the singular set. This must be considered a standard fact. Knowing the Green's function in our problem is essentially the same as knowing the Green's function for the Dirichlet problem for the ordinary Laplace equation on the product of a wedge of angle $2\pi \beta$ in $\bR^{2}$ with $\bR^{m-2}$ and, at least when $m=3$, this is a topic with a large classical literature (see for example \cite{kn:Carslaw}).  Equally, such polyhomogeneous expansions are prominent in the general theory of edge operators, as applied in \cite{kn:Mazzeo}. But, lacking an elementary reference for exactly the result we want, we will include a proof here. The proof involves traditional methods of separation of variables and a check on convergence.

We pause for a moment to recall some facts about Bessel functions. Our main reference is \cite{kn:WW}.  We fix $\nu\geq 0$. The Bessel equation for $f(z)$ is
\begin{equation}  f'' + z^{-1} f'+ (1-\nu^{2} z^{-2}) f =0. \end{equation}

 The Bessel function $J_{\nu}(z)$ is defined by a series expansion
\begin{equation}   J_{\nu}(z)=  \sum_{j=0}^{\infty}\frac{ (-1)^{j} (z/2)^{\nu+ 2j}}{ j! (\nu + j)!},\end{equation}
and satisfies the Bessel equation. (Here and below we use the notation $a!=\Gamma(a+1)$ for the generalised factorial function). The asymptotic behaviour for large real $z$ is $J_{\nu} \sim \sqrt{\frac{2}{\pi z}} \cos z$.    We define $J_{-\nu}$ by the same formula with $\nu$ replaced by $-\nu$. Then $J_{\nu}, J_{-\nu}$ are two solutions of the Bessel equation. They can be seen as roughly analogous to $\cos z, \sin z$. The linear combination
$$   h_{\nu}(z)= e^{\nu \pi i/2} J_{-\nu}(z)- e^{-\nu \pi i/2} J_{\nu}(z), $$
has the property that it decays rapidly at infinity on the upper half-plane: it is  roughly analogous to $e^{iz}$. 
We write $I_{\nu}(z)= e^{-\nu \pi i/2} J_{\nu}(iz)$ and
\begin{equation} K_{\nu}(z)= \frac{\pi}{2\sin(\nu \pi)} h_{\nu}(iz)= \frac{\pi}{2 \sin \nu \pi} \left( I_{-\nu}(z)-I_{\nu}(z))\right).\end{equation}
This formula can be extended to the case when $\nu$ is an integer by  taking a suitable limit. 

From (4), we have a convergent expansion
\begin{equation}I_{\nu}(z)= \sum_{j=0}^{\infty} \frac{1}{j!(\nu+j)!} \left( \frac{z}{2}\right)^{\nu+2j}, \end{equation}
and $I_{\nu}$
 has asymptotic behaviour for large positive $z$
\begin{equation}  I_{\nu}(z) \sim \frac{e^{z}}{\sqrt{2\pi z}}.  \end{equation}
The function  $K_{\nu}$ has the asymptotic behaviour for large positive $z$
\begin{equation}  K_{\nu}(z)\sim\sqrt{ \frac{\pi}{2z}} e^{-z}\ \end{equation}
but is unbounded near $z=0$. For $\nu>0$ 
$$  K_{\nu}(z)\sim \frac{(\nu-1)!}{2} \left( \frac{z}{2}\right)^{-\nu} \ \ z\rightarrow 0, $$
and $K_{0}(z)\sim -\log z$.
Our main tool will be the integral representation, 
\begin{equation}  K_{\nu}(z)= \frac{1}{2} \int_{-\infty}^{\infty} e^{ -z\cosh u + \nu u} du. \end{equation}

With this background in place, we proceed to analyse the Green's function by separation of variables. To begin with we argue formally, but in the end when we check convergence it will be clear that everything is watertight. Note that on grounds of symmetry we can write
$$ G(r,\theta, s; r',\theta',s')= \sum_{k\geq 0} G_{k}(r,r',R)\ \cos k(\theta-\theta'), $$
where $R=\vert s-s'\vert$. We want to find formulae for the $G_{k}$ and we  will usually write $\nu=c k$.
Our Laplace operator can be written as
$  \Delta_{g} \phi = \Delta_{\beta} \phi + \Delta_{\bR^{m-2}} \phi $ where
$\Delta_{\bR^{m-2}}$ is the ordinary Laplacian on $\bR^{m-2}$ and $\Delta_{\beta}$ is the operator in the plane defined by
$$\Delta_{\beta} \phi= \phi_{rr} + \frac{1}{r} \phi_{r} + \frac{1}{\beta^{2} r^{2}} \phi_{\theta \theta} . $$

This means that $\phi= J_{\nu} (\lambda r) e^{ik\theta}$ is an eigenfunction for $\Delta_{\beta}$, with  $\Delta_{\beta} \phi = -\lambda^{2}\phi$.  The Fourier-Bessel representation of a general function in terms of these eigenfunctions leads to a formula for the heat kernel associated to the operator $\Delta_{\beta}$ as 
$$   \sum_{k=0}^{\infty}  H_{k} \cos k (\theta-\theta') $$ where
$$  H_{k}(r,r')=\pi^{-1}  \int_{0}^{\infty} e^{-\lambda^{2} t} J_{\nu}(\lambda r) J_{\nu}(\lambda r')d\lambda.  $$
Now the heat kernel on a product is the product of the heat kernels, so the heat kernel of the Laplacian $\Delta_{g}$ on $\bR^{m}$ is
$$   (2\pi t)^{1-m/2} e^{-R^{2}/4t}\left( \sum H_{k}(r,r') \cos k(\theta-\theta')\right). $$
We  assume that $m\geq 3$. Then the Green's function can be obtained by integrating the heat kernel with respect to the time parameter. Thus

$$  G_{k}(r,r',R) =  \int_{0}^{\infty}\int_{0}^{\infty} (2\pi t)^{1-m/2}
e^{-\lambda^{2} t -R^{2}/4t} J_{\nu}(\lambda r) J_{\nu}(\lambda r')\ d\lambda dt. $$
Changing variable by $t=\frac{R}{2\lambda} e^{u}$ and using (9) we see that
$$  G_{k}= \frac{1}{(2\pi)^{m}} R ^{2-m/2}g_{k}$$ where
\begin{equation}  g_{k}= 2\int_{0}^{\infty} \lambda^{m/2-2} K_{m/2-2}( R \lambda  ) J_{\nu}(r\lambda ) J_{\nu}(r' \lambda ) d\lambda.
\end{equation}
The integral is convergent for all $r,r'$ provided that $R>0$. We get another representation by rotating the integration path. Suppose that $r<r'$ and write
$$  \sin (\nu \pi)  J_{\nu} (r' \lambda )= {\rm Im} (e^{-\nu \pi i/2} h_{\nu}(r' \lambda )). $$
Thus
$$  g_{k}= {\rm Im}\left[ \int_{0}^{\infty}\lambda^{m/2-2} K_{m/2-2}(2R \lambda) \frac{e^{-\nu \pi i/2}}{\sin \nu \pi} h_{\nu}(r'\lambda ) J_{\nu}(r\lambda) d\lambda \right]. $$
Because of the rapid decay of $h_{\nu}$ over the upper half plane we can rotate the integration path to the positive imaginary axis, which is the same as replacing $\lambda$ by $i\lambda$ in the integral. We get another expression
\begin{equation} g_{k}= 2 \int_{0}^{\infty} \lambda^{m/2-2} J_{m/2-2}(R\lambda) K_{\nu}(r'\lambda) I_{\nu}(r\lambda) d\lambda. \end{equation}
This integral converges for any $R$, provided that $r<r'$.

 \
 
We will now derive polyhomogeneous expansions for the Green's function in appropriate regions.  We need two elementary lemmas.

\begin{lem}
For $p,q\geq 0$ we have
$$  \int_{0}^{\infty} K_{p}(2x) x^{p+q} dx \leq (p+q)!. $$
\end{lem}

To see this, use the integral formula (9) to write the integral as
$$  I= \frac{1}{2} \int_{0}^{\infty}\int_{-\infty}^{\infty} e^{-2 x \cosh u + p u} x^{p+q} dx du. $$
Now change the order of integration and perform the $x$ integral to get
$$  I= \frac{(p+q)!}{2} \int_{-\infty}^{\infty} \frac{e^{pu}}{(2 \cosh u)^{p+q+1}} du. $$
Divide the integral into the two ranges $u\leq 0, u\geq 0$ and use the inequality
$2\cosh u\geq e^{-u}$ on the first and $2 \cosh u \geq e^{u}$ on the second. We get
$$  I \leq \frac{(p+q)!}{2} \left(\int_{-\infty}^{0} e^{(2p+q+1) u} du + \int_{0}^{\infty} e^{-(q+1) u} du \right). $$ The right hand side is
$$\frac{(p+q)!}{2} \left( \frac{1}{q+1}+ \frac{1}{2p+q+1}\right) \leq (p+q)!. $$

\begin{lem}
There is a universal constant $C$ such that for all $p,q\geq 1$ we have
 $$ \frac{(p+q)!}{p! q!} \leq C 2^{p+q}. $$
\end{lem}
When $p,q$ are integers this follows immediately from the binomial theorem, with $C=1$. The same argument applies when one of $p,q$ is an integer. Very likely we can always take $C=1$ but the author has not found this in texts so we give an ad hoc argument. Set $f(x)=  (1+ x)^{p} (1-x)^{q}$. We have an identity 
$$\int_{-1}^{1} f(x) dx = \frac{2^{p+q+1}}{p+q+1} \frac{p! q!}{(p+q)!}. $$
(See \cite{kn:WW} page 225.)
Since $(1+h^{-1})^{h}$ converges as $h\rightarrow \infty$ there is a universal constant $\delta>0$ such that if $0\leq x\leq 1/2q$ we have $(1-x)^{q} \geq \delta$, and if $-1/2p \leq x\leq 0$ we have $(1+x)^{q} \geq \delta$. So if $-1/2p\leq x\leq 1/2q$ we have $f(x) \geq \delta$. Thus the integral of $f(x)$ is at least
$\frac{\delta}{2}(p^{-1} + q^{-1})$. This gives
$$  2^{p+q} \frac{p! q!}{(p+q)!} \geq \frac{\delta}{4} (p+q+1) (p^{-1} + q^{-1}) \geq \delta/2. $$

First consider the representation arising from (11), when $r<r'$. It is convenient to normalise to $r'=2$. Write $J_{m/2-2}(x)= x^{m/2-2} F(x)$, so $F$ is bounded for positive real $x$. Then using the series representation (6) for $I_{\nu}$ and integrating term-by-term we get
\begin{equation}  G= \sum_{j,k} a_{j,k}(R) r^{\nu+2j} \cos k(\theta-\theta')\end{equation}
where
$$  a_{j,k}(R)= \frac{1}{2^{\nu+2j}} \frac{1}{j!(\nu+j)!} \int_{0}^{\infty} \lambda^{\nu+2j+m-4} F(R\lambda) K_{\nu}(2\lambda) d\lambda. $$
 Thus, by Lemma 1, 
 $$ \vert a_{j,k}(R)\vert \leq C' \frac{1}{2^{\nu+2j}} \frac{(\nu+2j+m-4)!}{j!(\nu+j)!}, $$
where $C'=\sup \vert F\vert$.  

Now, using Lemma 2,
$$  \vert a_{j,k}(R)\vert  \leq C C' \frac{(\nu+2j+m-4)!}{(\nu+2j)!} . $$
It is then elementary that the sum on the right hand side of (12) does converge absolutely provided that $r<1$.
For general $r'$ we can use the scaling behaviour to deduce that
$G= \left(\frac{2}{r'}\right)^{m-3} \sum a_{j,k}(\frac{2R}{r'}) \left(\frac{r}{2r'}\right)^{\nu+2j}\cos k(\theta-\theta')$ and the sum converges absolutely if $r<r'/2$.

For the other representation (10), we normalise to $R=1$. Then we expand both the Bessel functions $J_{\nu}(r\lambda ), J_{\nu}(r' \lambda )$ in powers of $\lambda$ and, arguing in a similar way, we have to consider a sum $\sum_{j,j',k} M_{j,j',k}$ where 
$$  M_{j,j',k} = \left(\frac{r}{2}\right)^{\nu+2j} \left(\frac{r'}{2}\right)^{\nu+2j'}  \frac{(p+2\nu+2j +2j')!}{j! j'! (\nu+j)!(\nu+j')!}. $$
 
 For $A,B,C,D>1$ we can write
 $$   \frac{(A+B+C+D)!}{A!B!C!D!}=\frac{(A+B+C+D)!}{(A+B)!(C+D)!}\frac{(A+B)!}{A!B!} \frac{(C+D)!}{C!D!}\leq C^{3} 2^{2(A+B+C+D)}, $$
by three applications of Lemma 2. Thus
$$   \frac{(2\nu+2j +2j')!}{j! j'! (\nu+j)!(\nu+j')!}\leq C^{3} 2^{2(2\nu+2j+2j')}. $$
Again, elementary arguments show that the sum of $M_{j,j',k}$ converges absolutely provided
$r,r'<1/2$. Scaling back: in the region $r,r'<R/2$ we get a convergent polyhomogeneous expansion
$$ G=  \sum b_{j,j',k}(R) r^{\nu+2j} (r')^{\nu+2j'} \cos k(\theta-\theta'). $$

We can now finish the proof of Proposition 4. We consider an open set $\Omega$ and $y$ not in $\Omega$. We know from standard elliptic regularity that we only need to verify the $\beta$-smoothness condition at  points $x_{0}$ in $\Omega\cap S$. Let $\Omega'$ be the ball of radius $d/10$ about $x_{0}$ where $d$ is the distance to the boundary of $\Omega$. We want to prove that  for any $y$ not in $\Omega$ the function $\Gamma_{y}$ restricted to $\Omega'$ is the composite of $\iota$ and a smooth function.    But for any such $y$ at least one of the two expansions above is valid over $\Omega'$. It is straightforward to see that the formal series we have considered define weak solutions of the equation characterising $\Gamma_{y}$ and since we have verified local convergence it follows that the sums represent valid formulae pointwise.  Now the series obviously define functions on appropriate neighbourhoods in $\bR^{m-2} \times \bR^{2}\times \bR$. That is, a polyhomogeneous series
$\sum a_{j,k}(s) r^{\nu+2j} \cos k(\theta)$ defines  a smooth function of $(s,\zeta, \rho)$ by 
$$ a_{j,k}(s)  \rho^{j} {\rm Re}( \zeta^{k}), $$
which gives the required factorisation when we set $\rho=r^{2}, \zeta= r^{c} e^{i\theta}$. Similar considerations show that $\Gamma_{y}$ varies continuously in the desired sense as $y$ varies.

\section{Application to K\"ahler-Einstein equations}

\subsection{Flat model}
Write $\bC_{\beta}$ for the Riemannian manifold with underlying space $\bC$, on which we take a standard co-ordinate $\zeta$,  and with the singular metric associated to the $2$-form  $ \beta^{2} \vert \zeta \vert^{2(\beta-1)}id\zeta d\overline{\zeta}$.
Then the map $\zeta= r^{c} e^{i\theta}$ (recall that we write $c=\beta^{-1}$) gives an isometry from the standard cone metric $dr^{2} + \beta^{2} r^{2} d\theta^{2}$ to $\bC_{\beta}$. Likewise, when $m=2n$ we get an isometry between the singular metric we considered above on $\bR^{2m}$ and the Riemannian product $\bC_{\beta}\times \bC^{n-1}$. Let $\sigma_{1},\dots, \sigma_{n-1}$ be standard complex co-ordinates on $\bC^{n-1}$. Thus we have two natural systems of co-ordinates $(r,\theta, \sigma_{a})$ and $(\zeta, \sigma_{a})$

We  consider the $ \dbd$-operator on the complement of $S$, mapping functions to (1,1) forms. Set $\epsilon = dr+ i  \beta r d\theta$. Then, up to a factor of $\sqrt{2}$, the forms  $\epsilon, d\sigma_{1}, \dots, d\sigma_{m-1}$ give an orthonormal basis for the $(1,0)$ forms at each point. We should keep in mind that $\epsilon$ is {\it not} a holomorphic $1$-form , although $c r^{c-1} e^{i\theta} \epsilon=d\zeta$ is. Now   take a trivialisation of the $(1,1)$ forms by sections
$$  d\sigma_{a} \wedge d\overline{\sigma}_{b}\ , \ d\sigma_{a}\wedge  \epsilon\ ,\  d\overline{\sigma}_{a} \wedge \overline{\epsilon}\ , \ \epsilon\wedge \overline{\epsilon}. $$
 Up to scale factor, this is a unitary trivialisation. With  respect to this trivialisation the components of $\dbd$ are all operators $D$ of the kind  considered above, except for 
$$D_{0}= i \left( r^{-1} \frac{\partial}{\partial r}(r \frac{\partial }{\partial r}) + \beta^{-2} r^{-1} \frac{\partial^{2}}{\partial\theta^{2}} \right). $$
which is the $i \epsilon\wedge \overline{\epsilon}$ component of $\dbd$.
Of course this is just the Laplacian  $\Delta_{\beta}$ in the $\bR^{2}$ variable, with respect to the singular metric. So
$$  \Delta_{g}= D_{0} + \Delta_{\bC^{n-1}}$$
in an obvious notation. 
Since, by definition, $\Delta_{g} G\rho=\rho$ we can write
$$  D_{0} G \rho= \rho- \Delta_{\bC^{n-1}} \rho. $$
The operator $\Delta_{\bC^{n-1}}$ is a sum of terms of the form allowed in Section 2 so  we  get  a H\"older estimate on
$D_{0}G \rho$ and hence on $\dbd G\rho$.  
So we have
\begin{cor} Suppose $\alpha<\mu=(\beta^{-1}-1)$. Then there is a constant $C$ depending only on $m,\beta,\alpha$ such that for all $\rho\in C^{\infty}_{c}$ we have
$$  [ \dbd (G\rho) ]_{\alpha} \leq C [\rho]_{\alpha}, $$
where the left hand side is interpreted using the trivialisation above. 
\end{cor}

\

Notice that it follows from our discussion of the Green's function that the components of $\dbd G\rho$ corresponding to the basis elements $\epsilon \wedge d\sigma_{a}$ tend to zero on the singular set $S$.

\subsection{Further local theory}
Corollary 1 expresses the essential fact that we are after, but for applications
we need a variety of other statements which will be set out here. One detail is that the smooth functions are not dense in H\"older spaces. But any $C^{,\alpha}$ function can be approximated by smooth functions in the norm of $C^{,\underline{\alpha}}$ for any  $\underline{\alpha}<\alpha$. So in the end this complication becomes irrelevant and we will ignore it. Suppose that $\rho$ is a $C^{,\alpha}$ function with support in the unit ball $B\subset \bC_{\beta}\times \bC^{n-1}$. Then $\dbd G\rho$ is $C^{,\alpha}$ and the same estimate  as in Corollary 1 holds. 
As in our discussion of $\dbd$, we say that the derivative of a function $f$ is in $C^{,\alpha}$ if the components $\frac{\partial f}{\partial r}, r^{-1} \frac{\partial f}{\partial \theta}$ and $\frac{\partial f}{\partial s_{i}}$ are $C^{,\alpha}$.  
Similar arguments  to those  of Section 2 (but easier)show that in the situation above $G\rho$ and $\nabla G\rho$ are in $C^{,\alpha}$. In fact the same argument show that $G\rho, \nabla G\rho$ are in $C^{,\overline{\alpha}}$ for any $\overline{\alpha}$ with $\overline{\alpha}<\mu$ and we have an estimate
\begin{equation}    [G\rho]_{\overline{\alpha}}+[\nabla G\rho]_{\overline{\alpha}}\leq C [\rho]_{\alpha}. \end{equation}
Taking $\overline{\alpha}>\alpha$ we get a compactness result: for a sequence  $\rho_{i}$, supported on $B$ and bounded in $C^{,\alpha}$, there is a subsequence $\{i'\}$ such that $G\rho_{i'}$ and $ \nabla G\rho_{i'}$ converge in $C^{,\alpha}$ over compact sets. 
Notice also that, as in the remark following Corollary 1, the components of $\nabla G\rho$ corresponding to the derivatives  $\frac{\partial f}{\partial r}, r^{-1} \frac{\partial f}{\partial \theta}$ tend to zero on the singular set.

Now consider the situation where we have a function $\phi\in C^{,\alpha}(B)$ such that $\Delta \phi$, defined pointwise outside $S$, is also $C^{,\alpha}$. Applying standard elliptic estimates in small balls in the complement of $S$ we see that $\vert \nabla \phi\vert =O(r^{-1+\alpha})$ near the singular set. One easy consequence is that $\Delta \phi$, defined pointwise as above, agrees with the weak, distributional, notion. For another we take a smooth cut-off function $\chi$ of compact support in $B$, equal to $1$ on some interior region $B'$ and  with $\Delta \chi$ smooth. Then $\Delta(\chi \phi)$ is in $L^{q}$ so $G\Delta (\chi \phi)$ is defined. It follows that
$\chi \phi = G\rho_{1}+G\rho_{2}$ where $\rho_{1}= (\Delta \chi) \phi+ \chi \Delta \phi$ and $\rho_{2}= 2 \nabla\chi.\nabla\phi$. Away from the support of $\nabla \chi$ it is clear that $G\rho_{2}$ is in $C^{,\alpha}$. Thus we see that $\dbd \phi$ is locally in $C^{,\alpha}$ and we obtain  interior estimates of the form
\begin{equation}   [\dbd \phi]_{\alpha, B'} \leq C \left( [\Delta \phi]_{\alpha, B}+[\phi]_{\alpha,B}\right), \end{equation}
where $C$ depends on $B'$. Similarly we get
\begin{equation} [\nabla \phi ]_{\overline{\alpha}, B'} \leq C \left( [\Delta \phi]_{\alpha, B}+[\phi]_{\alpha,B}\right)\end{equation}
for any $\overline{\alpha}<\mu$. 

Now let $\eta$ be a $C^{,\alpha}$ section of the bundle of $(1,1)$ forms, in the sense we have defined, and consider the operator $\Delta_{\eta}\phi =\Delta \phi + \eta. \dbd \phi$. We suppose first that $\eta$ is supported on $B$ and is sufficiently small in $C^{,\alpha}$. It follows from the usual Neumann series argument that we can invert $\Delta_{\eta}$  and that an estimate corresponding to Corollary 1 holds. Then we can extend all the results above to $\Delta_{\eta}$. As usual, if we have any $\eta$ which vanishes at the origin we can reduce to the situation where $\eta$ is small and of compact support by dilation and multiplying by a cut-off function and thus obtain the interior estimate near  the origin. 

\subsection{Global set-up}

Let $X$ be a compact K\"ahler manifold and $D\subset X$ be a smooth hypersurface.
Let $\Lambda\rightarrow X$ be the holomorphic line bundle associated to $D$, so there is a  section $s$ of $\Lambda$ cutting out $D$. Let $h_{\Lambda}$ be any smooth hermitian metric on $\Lambda$ and write
$$  \chi= i  \partial \db \vert s \vert^{2\beta}_{h_{\Lambda}}. $$
 Let $\Omega_{0}$ be a smooth K\"ahler metric on $X$. Then we have
\begin{lem}
For sufficiently small $\delta>0$ the $(1,1)$ form $\omega_{0}=\Omega_{0}+\delta \chi$ is positive on $X\setminus D$. The metric we obtain is independent of the choices of $\Omega_{0}, h_{\Lambda}, \delta$ up to quasi-isometry.
\end{lem}
This is elementary to check and we omit the proof. If we choose standard complex co-ordinates  $\zeta, \sigma_{a}$ around a point of $D$, so that $D$ is defined by the equation $\zeta=0$, 
then $\vert s\vert^{2}_{h_{\Lambda}}= F \vert \zeta \vert^{2}$ where $F$ is a smooth positive function of $\zeta, \sigma_{a}$. Thus
\begin{equation} \chi =  (\dbd F^{\beta}) \vert \zeta \vert^{2\beta} +
i \beta \vert \zeta \vert^{2(\beta-1)} \left( \zeta \partial F^{\beta} d \ozeta - \ozeta \db F^{\beta} d\zeta\right) + \beta^{2} F^{\beta} \vert \zeta \vert^{2(\beta-1)}  i d\zeta d\ozeta. \end{equation}

Lemma 3 implies that there is a well-defined notion of a H\"older continuous function, with exponent $\alpha$, on $X\setminus D$, using the singular metric. If we take a standard local complex co-ordinate system $\zeta, \sigma_{a}$ as above and then set $z = \zeta \vert \zeta\vert^{\beta-1}$  then this becomes the ordinary notion of H\"older continuity in terms of the co-ordinates $z,\sigma_{a}$.
 We write $C^{,\alpha,\beta}$, or sometimes just $C^{,\alpha}$ for these functions on $X$. Now we want to go on to define H\"older continuous differential forms. With a fixed metric $h_{\lambda}$ as above, define the $(1,0)$-form 
$$ \eta = \partial \vert s \vert^{\beta}. $$
Then we say that a $(1,0)$ form on $X\setminus D$ is H\"older continuous with exponent $\alpha$ if and only if it can be written as
  $$    f_{0} \eta + f_{1} \pi_{1} + \dots +f_{N} \pi_{N}, $$
  where $f_{i}\in C^{,\alpha,\beta}$ and $\pi_{i}$ are smooth forms (in the ordinary sense) on $X$. 

\begin{lem} If $\alpha<\mu$ this notion is independent of the choice of metric $h_{\Lambda}$. 
\end{lem}
If we make another choice of metric we get another form $\eta'= \partial( f \vert s\vert^{\beta})$ for a smooth positive function $f$. Then
$$\eta'= \partial f \vert s\vert^{\beta} + f \eta.$$
The function $\vert s\vert^{\beta}$ is in $C^{,\alpha}$ so $\eta '$ can be written in the stated form, and the result follows immediately. Similarly ones sees that, in standard co-ordinates $z=re^{i\theta}$, the H\"older continuous (1,0)- forms are just those of the shape $f_{0}\epsilon + \sum_{a=1}^{n-1} f_{a} d\sigma_{a}$, where $\epsilon= dr+ i\beta r d\theta$, as in the previous subsection,   the co-efficients $f_{0}, f_{1} \dots$ are in $C^{,\alpha}$ and $f_{0}$ vanishes on the singular set.
Similarly, we can give a global definition of a space of  H\"older continuous $(1,1)$ forms which reduces in local co-ordinates $(r,\theta,\sigma_{a})$ to those of the shape
\begin{equation}    m i  \epsilon \overline{\epsilon} + \sum m_{ab} d\sigma_{a}d\overline{\sigma}_{b} + m_{a} \epsilon d\overline{\sigma}_{a} + \overline{m}_{a} \overline{\epsilon}d\sigma_{a}
\end{equation}
where $m, m_{ab}, m_{a}$ are $C^{,\alpha}$ and the $m_{a}$ vanish on the singular set.   

Now we define ${\cal C}^{2,\alpha,\beta}$ to be the space of (real-valued) functions $f$ on $X\setminus D$ with $f, \partial f, \dbd f $ all H\"older continuous with exponent $\alpha$. This is the analogue of the usual H\"older space $C^{2,\alpha}$ but there is an important difference that we are not asserting that {\it all} second derivatives of $f$ are in $C^{,\alpha}$. We can define norms on $C^{,\alpha,\beta}, {\cal C}^{2,\alpha,\beta}$ in the usual way, making them Banach spaces.   If $\omega_{0}$ is a singular metric, as constructed above, we have a space of {\it H\"older continuous K\"ahler metrics} of the form $\omega_{\phi}= \omega_{0}+\dbd \phi$ where $\phi\in {\cal C}^{2,\alpha,\beta}$ and we require that $\omega_{\phi}\geq \kappa \omega_{0}$ on $X\setminus D$, for some $\kappa>0$. It is easy to check that this space of metrics is independent of the choice of $\omega_{0}$.

Let $\omega$ be a H\"older continuous K\"ahler metric as above. In local co-ordinates the metric is described by co-efficients as in (17).  All of these have limits along the singular set and by definition the limits of the $m_{a}$ are zero. The limits of the $m_{ab}$ obviously define a $C^{,\alpha}$ metric on $D$ and the limit of the power $m^{1/\beta}$ is intrinsically a $C^{,\alpha}$ Hermitian metric on the restriction of the line bundle $\Lambda$ to $D$ (which is identified with the normal bundle of $D$ in $X$). Given any point $p\in D$ it is clear that we can choose a standard co-ordinate system centred at $p$ so that the $m=1$ and $m_{\alpha \beta}=\delta_{\alpha \beta}$ at this point.
Now write $\Delta$ for the Laplace operator of the metric $\omega$. Since it is given by an algebraic contraction of $\dbd$ it appears, in these local co-ordinates, in the form $\Delta_{\eta}$ considered in the previous subsection, and $\eta$ vanishes at $p$. So we can apply the results there to obtain interior estimates and inversion operators in sufficiently small balls about this point. From here we can carry through the usual arguments to   obtain a parametrix for $\Delta$ over all of $X$. In this way we obtain
\begin{prop}
If $\alpha<\mu=(\beta^{-1}-1)$   the inclusion ${\cal C}^{2,\alpha,\beta}\rightarrow C^{,\alpha,\beta}$ is compact. If $\omega$ is a ${\cal C}^{,\alpha,\beta}$ K\"ahler metric on $(X,D)$ then the Laplacian of $\omega$ defines a Fredholm map $\Delta:{\cal C}^{2,\alpha,\beta}\rightarrow C^{,\alpha,\beta}$. 
\end{prop}

From now on we restrict attention to the case when $X$ is a Fano manifold,  $[\Omega_{0}]= 2\pi c_{1}(X)$ and $D$ is in the linear system $-K_{X}$. We can regard $\Omega_{0}$ as the curvature form associated to a smooth metric  on the dual of $K_{X}$. Then $\omega$ is the curvature form of  a singular metric $h_{0}$ on this line bundle and  any $\phi\in {\cal C}^{2,\alpha,\beta}$ defines another metric $\vert \ \vert_{\phi}=e^{\phi} h_{0}$.
  We identify $\Lambda_{D}$ with $K_{X}^{-1} $, so $s$ is a section of $K_{X}^{-1}$.
If $\omega_{\phi}$ is any K\"ahler metric on $X\setminus D$ its Riemannian volume form can be regarded as an element of $K_{X}\otimes \overline{K}_{X}$ so we get  a function
$   s\otimes \overline{s} \Vol_{\omega}$ on $X\setminus D$.  We say the metric is K\"ahler-Einstein if
\begin{equation}    s\otimes \overline{s} \Vol_{\omega}= \vert s\vert_{\phi}^{2\beta}. \end{equation}
 If this holds then, by standard elliptic regularity, $\phi$ is smooth on $X\setminus D$ and satisfies $\Ric(\omega_{\phi})= \beta \omega_{\phi}$.

We expect there to be a detailed regularity theory for these Kahler-Einstein metrics around the singular divisor, as outlined by Mazzeo in \cite{kn:Mazzeo}. We will  leave most of the discussion of this to another paper but we want to observe here that the metrics are \lq\lq smooth in tangential directions''.
In a local co-ordinate system $(z,\sigma_{a})$ we can choose a local K\"ahler potential $\psi$ so that the equation becomes
$$   (\dbd \psi)^{n}= e^{\beta\psi}. $$
Let $\psi'$ be a derivative with respect to the real or imaginary part of any $\sigma_{a}$.  Then $\psi'$ satisfies a linear equation $(\Delta + \beta)\psi'=0$, so it follows that $\dbd \psi'$ is $C^{,\alpha}$. Repeating the argument, we find that all multiple derivatives in these directions satisfy this condition. In particular, the induced metric on $D$ and the metric induced on the restriction of $\Lambda$ to $D$ are both smooth. 

Another simple fact is that a solution of our K\"ahler-Einstein equation which is in ${\cal C}^{2,\alpha,\beta}$ for some $\alpha>0$ lies in ${\cal C}^{2,\overline{\alpha},\beta}$ for all $\overline{\alpha}<\mu=\beta^{-1}-1$: thus the theory is independent of the choice of exponent $\alpha$.

\subsection{Deforming the cone angle}

For a Fano manifold $X$ and smooth $D\in \vert -K_{X}\vert$ as above we have:
\begin{thm} 
Let $\beta_{0} \in (0,1), \alpha< \mu_{0}=\beta_{0}^{-1}-1$ and suppose there is a ${\cal C}^{2,\alpha, \beta_{0}}$ solution $\omega$ to the K\"ahler-Einstein equation (18) on $(X,D)$, with $\beta=\beta_{0}$. If there are no nonzero holomorphic vector fields on $X$ which are tangent to $D$ then for $\beta$ sufficiently close to $\beta_{0}$ there is a ${\cal C}^{2,\alpha,\beta}$ solution to (18) for this cone angle.
\end{thm}

It seems likely that the condition on holomorphic vector fields is always satisfied, by general results from algebraic geometry, but the author has not gone into this. In any case it is not a serious restriction. 

The proof of the theorem follows standard general lines. Having set up a linear theory, we can deform the solutions to the nonlinear equation using an implicit function theorem, provided that the linearised operator is invertible. However there are some complications, for example due to the fact that the function spaces depend on $\beta$. We have seen that the solution $\omega$ defines a smooth
metric on $\Lambda$ over $D$. We extend this to a smooth metric, which we will write as $\Vert\ \Vert$, on $\Lambda$ over $X$. This is not to be confused with the singular metric, which we will write as $\vert \ \vert$, whose curvature is $\omega$. Now for $\beta$ near to $\beta_{0}$ we define
$$  \omega_{\beta}= \omega + \dbd(\Vert s\Vert^{\beta}- \Vert s \Vert^{\beta_{0}}), $$
so $\omega_{\beta_{0}}=\omega$. In other words, $\omega_{\beta}$ is the curvature form of the singular metric on $K_{X}^{-1}$ with
$$ \vert s\vert_{\beta}^{2}= \exp( \Vert s\Vert^{\beta}-\Vert s\Vert^{\beta_{0}}) \vert s \vert^{2}. $$
 Set $$ k_{\beta}=  \vert s \vert_{\beta}^{-2\beta} s\otimes \overline{s} \ \Vol_{\omega_{\beta}} . $$
Thus $k_{\beta_{0}}=1$, since $\omega$ solves the K\"ahler-Einstein equation. We state three Propositions.
\begin{prop}  $$\Vert k_{\beta}-1 \Vert_{C^{,\alpha,\beta}}\rightarrow 0$$ as $\beta\rightarrow \beta_{0}$.
\end{prop}
Write $\Delta_{\beta}$ for the Laplace operator of $\omega_{\beta}$.
\begin{prop}
If $\Delta_{\beta_{0}}+\beta_{0}: {\cal C}^{2,\alpha,\beta_{0}}\rightarrow C^{,\alpha,\beta_{0}}$ is invertible then for $\beta$ close to $\beta_{0}$ the operator $\Delta_{\beta}+\beta:{\cal C}^{2,\alpha,\beta}\rightarrow C^{,\alpha,\beta}$ is also invertible and the operator norm of its inverse is bounded by a fixed constant independent of $\beta$.
\end{prop}

The statements of  Propositions 6 and 7 are not completely precise. There are many ways of defining norms on $C^{\,\alpha,\beta}, {\cal C}^{2,\alpha,\beta}$, all of which are equivalent for fixed $\beta$. But what we need here is a definite family of norms, for example defined using a fixed system of co-ordinate charts. But we hope that the details of such a definition will be clear to the reader and do not need to be spelled out.

\begin{prop}If $\Delta_{\beta_{0}} +\beta_{0}$ is not invertible then there is a non-trivial holomorphic vector field on $X$ tangent to $D$.
\end{prop}

Given these three results, the proof of Theorem 2 is a standard application of the implicit function theorem. 

\

We begin with the proof of Proposition 6. This is completely elementary, but the set-up is a little complicated. As a first simplification we reduce to considering convergence with respect to the H\"older norm defined by the fixed parameter $\beta_{0}$. That is to say, for any $\beta$ we are considering a standard chart $\chi_{\beta}$ mapping a neighbourhood of $0$ in $\bC\times \bC^{n-1}$ to $X$ and the functions in $C^{,\alpha,\beta}$ are those which pull back by $\chi_{\beta}$ to ordinary $C^{,\alpha}$ functions. The composite $\eta_{\beta, \beta_{0}}=\chi_{\beta}^{-1} \circ \chi_{\beta}$ is the map defined by $(re^{i\theta}, s)\mapsto (r^{\lambda} e^{i\theta}, s)$, where $\lambda=\beta_{0}/\beta$. If $\beta>\beta_{0}$ this is not Lipschitz so the notions of H\"older continuity are different. However, $\eta_{\beta,\beta_{0}}$ is $\beta_{0}/\beta$-H\"older and this means that it pulls $C^{,\alpha}$ functions back to $C^{,\alpha \beta_{0}/\beta}$ functions. Since we are always free to adjust $\alpha$ a little and since we can take $\beta_{0}/\beta$ arbitrarily close to $1$, we see that it suffices to prove that  

$$\Vert k_{\beta}-1 \Vert_{C^{,\alpha,\beta_{0}}}\rightarrow 0$$ as $\beta\rightarrow \beta_{0}$.

We will use another, similar, elementary observation  below.  Supppose that $f_{i}$ is a sequence of functions on the ball in $\bC\times \bC^{n-1}$ converging to a limit $f_{\infty}$ in $C^{,\alpha}$ and with $f_{i}$ all vanishing on the singular set $\{r=0\}$. Suppose that $0<\epsilon_{i}\leq \epsilon<\alpha$ and $\epsilon_{i}\rightarrow 0$ as $i\rightarrow \infty$. Then the functions $r^{-\epsilon_{i}} f_{i}$ are H\"older with exponent $\alpha-\epsilon$ and converge in this sense to $f_{\infty}$ as $i\rightarrow \infty$. 

\

With these remarks in place we can begin the proof. We work in a standard local co-ordinate system $\zeta, \sigma_{a}$ chosen so that section $s$ is given by $$ s= \zeta (d\zeta d\sigma_{1} \dots d\sigma_{n})^{-1}. $$
Then $\Vert s\Vert^{2} = F \vert \zeta\vert^{2}$, where $F$ is smooth strictly positive function of $\zeta, \sigma_{a}$. Now write, as in (16), 
$$   \dbd (F^{\beta} \vert \zeta\vert^{2\beta})= F^{\beta} \vert \zeta \vert^{2\beta-2} \tau + V_{\beta}, $$
say, where $\tau= i d\zeta d\ozeta$. Of course we can write down a formula for $V_{\beta}$, although it is a little complicated. The point to emphasise is that this just depends on the smooth function $F$ and $\beta$. All we need to know is that the (1,1)-forms $V_{\beta}$ are $C^{,\alpha,\beta_{0}}$ forms for $\beta$ close to $\beta_{0}$, they all vanish on the singular set and they converge to $V_{\beta_{0}}$ in this H\"older space sense as $\beta\rightarrow \beta_{0}$. We leave the reader to verify these assertions by straightforward calculation.

Now we can write
$$  \omega_{\beta}= F^{\beta} \vert \zeta \vert^{2\beta-2} \tau + V_{\beta}+\Omega, $$ where $\Omega$ is independent of $\beta$. Thus in our standard co-ordinates $r, \theta, \sigma_{a}$ the form $\Omega$ has $C^{,\alpha}$ co-efficients and all co-efficients tend to zero on the singular set except those involving $d\sigma_{a}d\overline{\sigma}_{b}$. 

Recall that $k_{\beta}= \vert s\vert_{\beta}^{-2\beta} s\otimes \overline{s} \Vol(\omega_{\beta})$. 
By the definition of our class of Holder continuous metrics we can write
$$\vert s \vert^{2}_{\beta_{0}}= \Vert s\Vert^{2} \exp(\psi+ \Vert s\Vert^{2\beta_{0}}), $$
where $\psi$ is ${\cal C}^{2,\alpha,\beta_{0}}$. From this we get
$$   \frac{k_{\beta}}{k_{\beta_{0}}} = \Vert s \Vert^{2(\beta_{0}-\beta)} \exp((\beta-\beta_{0}) \psi +\beta\Vert s\Vert^{2\beta}-\beta_{0} \Vert s\Vert^{2\beta_{0}}) \frac{\Vol(\omega_{\beta})}{\Vol(\omega_{\beta_{0}})}. $$
Writing $\Vert s\Vert^{2} = F \vert \zeta \vert^{2}$ we get
$$   \frac{k_{\beta}}{k_{\beta_{0}}}= \vert \zeta \vert^{2(\beta_{0}-\beta)} H_{\beta} \frac{\Vol(\omega_{\beta})}{\Vol(\omega_{\beta_{0}})}, $$
where $H_{\beta}$ tends to $1$ in $C^{,\alpha,\beta_{0}}$ as $\beta\rightarrow \beta_{0}$. So it suffices to prove that
 $\vert \zeta\vert^{2(\beta_{0}-\beta)}\frac{\Vol(\omega_{\beta})}{\Vol(\omega_{\beta_{0}})}$ also tends to $1$ in this sense. 

For simplicity, to explain the argument, let us suppose that $n=2$. Write $V_{\beta}+\Omega=\Omega_{\beta}$. So $\Omega_{\beta}$ are $(1,1)$ forms which vary continuously in  $C^{,\alpha,\beta_{0}}$ for $\beta$ close to $\beta_{0}$.
We take the standard volume form $J_{0}$ in our co-ordinates $(r,\theta,\sigma_{0})$ to be $J_{0}= \tau\wedge d\sigma_{1} \dots d\sigma_{n-1}d\osigma_{1}\dots d\osigma_{n-1}$. Then $\omega_{\beta_{0}}^{2}/J_{0}\geq \kappa>1$. Now since $\tau^{2}=0$ we have
$$\omega_{\beta}^{2}= F^{\beta} \vert \zeta\vert^{2(\beta-\beta_{0})} \tau \wedge\Omega_{\beta} + \Omega_{\beta}^{2}, $$
so, writing $r=\vert \zeta\vert^{\beta}$, 
$$  \vert \zeta\vert^{2(\beta_{0}-\beta)} \omega_{\beta}^{2}= F^{\beta} \tau\wedge \Omega_{\beta} +r^{2(1-\beta/\beta_{0})} \Omega_{\beta}^{2}. $$
The crucial thing is that $\Omega_{\beta}^{2}/J_{0}$ vanishes on the singular set. Thus we can apply the observation about multiplication above to see that, after slightly adjusting $\alpha$, the product $r^{2(1-\beta/\beta_{0})} \Omega_{\beta}^{2}/J_{0}$ converges in $C^{,\alpha,\beta_{0}}$ as $\beta$ tends to $\beta_{0}$. This completes the proof of Proposition 6.
\

Next we consider Proposition 8. The first step is to establish a Fredholm alternative: if $(\Delta_{\beta_{0}}+\beta_{0})$ has no kernel in ${\cal C}^{2,\alpha,\beta_{0}}$ then it is surjective (i.e. the Fredholm index is zero). For the corresponding $L^{2}$ theory this is straightforward, so what one needs to know is that if $\rho$ is in $C^{,\alpha}$ and $f$ is a weak solution of the equation
$(\Delta+\beta_{0})f=\rho$ then $f$ is in ${\cal C}^{2,\alpha,\beta_{0}}$. By the results of (4.2), this will be true if we can show that $f$ is in $C^{,\alpha,\beta_{0}}$ and this follows from the general theory developed in \cite{kn:GT}, Chapter 8. Granted this, the proof of Proposition 8 comes down to showing that a non-trivial solution of $(\Delta_{\beta_{0}}+\beta_{0})f=0$ defines a non-trivial holomorphic vector field on $X$, tangent to $D$. Of course this is standard material in the ordinary, non-singular, case.To simplify notation write $\Delta_{\beta_{0}}=\Delta$.   We write $\cD$ for the operator $\db\circ \grad$ over $X\setminus D$, where $\grad f$ denotes the gradient vector field of $f$ with respect to the metric and $\db$ is the $\db$-operator on vector fields. The fact that the Ricci curvature of $\omega_{\beta_{0}}$ is $\beta_{0}\omega_{\beta_{0}}$ gives an identity
\begin{equation}   \db^{*} \cD f= \grad(\Delta f+\beta_{0} f). \end{equation}
For $\epsilon>0$, let $X_{\epsilon}$ be the complement of a tubular neighbourhood of $D$ in $X$, modelled in standard local co-ordinates on the region $\{r\geq \epsilon\}$. Suppose that $(\Delta +\beta_{0})f=0$. We take the inner product of (19) with $\grad f$ and integrate by parts to get
\begin{equation}  \int_{X_{\epsilon}}\vert \cD f \vert^{2} = \int_{\partial X_{\epsilon}} \cD f * \grad f, \end{equation}
where $*$ denotes a certain bilinear algebraic operation. What we need to see is that the boundary term tends to $0$ with $\epsilon$. From that we see that $\grad f$ is a holomorphic vector field on $X\setminus D$. We know that the radial derivative $\frac{\partial f}{\partial r}$ is $O(r^{\alpha})$ for and this translates into the fact that the $\frac{\partial}{\partial \zeta}$ component of the vector field, in holomorphic co-ordinates, is $O(\vert \zeta\vert^{\alpha \beta-\beta+1})$. The $\frac{\partial}{\partial \sigma_{a}}$ components are bounded. If $\alpha$ is sufficiently close to $\mu=\beta^{-1}-1$ then $\alpha \beta-\beta+1$ is positive and this implies that the  vector field extends holomorphically across $D$ and is tangent to $D$.

So the real task is to check that the boundary term in (20)   tends to zero with $\epsilon$. For this we use
\begin{lem}
With the notation above, $\vert \cD f\vert =O(r^{\alpha-1})$, where $r$ is the distance (in the metric $\omega_{\beta_{0}}$) to the divisor $D$.
\end{lem}
Assuming this Lemma it follows that the integrand $\cD f *\grad f$ is $O(r^{\alpha})$, because $\grad f$ is bounded so the boundary integral is $O(r^{\alpha})$ and the volume of $\partial X_{\epsilon}$ is $O(r)$. 

To prove the Lemma we can work in a local chart and there is no loss in taking $r$ to be the radial co-ordinate as before. Given a point $p$ with radial co-ordinate $r_{0}$ we consider a small ball $B_{0}$ of radius $ h r_{0}$ centred at $p$ on which we can identify the model cone metric with the flat metric (so $h$ is a fixed small number depending on $\beta_{0})$. We re-scale this small ball to a unit ball $B\subset \bC^{n}$. The K\"ahler-Einstein metric $\omega_{\beta_{0}}$ re-scales to a K\"ahler-Einstein metric $\tilde{\omega}$ on $B$. The fact that $\omega_{\beta_{0}}$ is $C^{,\alpha}$ means that the $C^{0}$ difference between $\tilde{\omega}$ and a Euclidean metric on $B$ is $O(r_{0}^{\alpha})$. Now standard elliptic regularity for the K\"ahler-Einstein equations implies that the derivative of $\tilde{\omega}$ is also $O(r_{0}^{\alpha})$ on an interior ball. Scaling back, we see that the derivative of $\omega_{\beta_{0}}$ is $O(r_{0}^{\alpha-1})$ at $p$. 

Now consider our function $f$ with $(\Delta +\beta_{0})f=0$. We know that the radial derivative of $f$ is $O(r^{\alpha})$ and the tangential derivatives are in $C^{,\alpha}$. Given $p$ as above, let $f_{0}$ be the $\bR$-linear function of the co-ordinates $\sigma_{a}$ defined by the tangential derivative of $f$ at $p$. Thus the derivative of $g=f-f_{0}$ is $O(r_{0}^{\alpha})$ over $B_{0}$ and the variation of $g$ over $B_{0}$ is $O(r_{0}^{\alpha+1})$. We also have $\Delta g=-\beta_{0} f $  since $\Delta f_{0}=0$. By the same kind of argument as before, rescaling and using standard elliptic estimates, we see that $\cD g$ is $O(r_{0}^{\alpha-1})$ at $p$. On the other hand $\cD f_{0}= \db(\grad f_{0})$ and the definition of $\grad f_{0}$ involves the metric tensor $\omega_{\beta_{0}}$. From this we see that $\vert \cD f_{0}\vert$ is bounded by a fixed multiple of the derivative of the metric tensor and so is $O(r_{0}^{\alpha-1})$ by the preceding discussion. Hence $\cD f= \cD g +\cD f_{0}$ is $O(r_{0}^{\alpha-1})$ as required.

Finally we turn to Proposition 8, but here we will be very brief since  nothing out of the ordinary is involved. By the Fredholm alternative, it suffices to show that if $\beta_{i}$ is a sequence converging to $\beta_{0}$ and if $f_{i} $ are functions with $\Vert f_{i}\Vert_{{\cal C}^{2,\alpha,\beta_{i}}}=1$ but $\Vert (\Delta_{\beta_{i}} +\beta_{i}) f_{i}\Vert_{C^{,\alpha,\beta_{i}}}\rightarrow 0$ as $i\rightarrow \infty$ then there is a nontrivial solution to the equation
$(\Delta_{\beta_{0}}+\beta_{0}) f=0$. To do this one applies elementary observations about the family of metrics $\omega_{\beta}$, like those in the proof of Proposition 6, and standard arguments to get uniform estimates, independent of $i$.

\section{Model Ricci-flat solutions}
\subsection{Digression in four-dimensional Riemannian geometry}
Suppose that we have six $2$-forms $\omega_{1},\omega_{2},\omega_{3}, \theta_{1}, \theta_{2}, \theta_{3} $ on a 4-manifold which satisfy the equations
\begin{equation} \omega_{i}\wedge \omega_{j}= V \delta_{ij}\ ,\ \theta_{i}\wedge \theta_{j}=-V \delta_{ij} \ ,\ \omega_{i}\wedge \theta_{j}= 0\end{equation}
where $V$ is a fixed volume form. There is a unique Riemannian metric such that the $\omega_{i}$ form an orthonormal basis for the self-dual forms $\Lambda^{+}$ and $\theta_{j}$ for the anti-self-dual forms $\Lambda^{-}$. We want to discuss the Levi-Civita connection of this metric, viewed as a pair of connections on the bundles $\Lambda^{+}, \Lambda^{-}$ (that is,  using the local isomorphism between $SO(4)$ and $SO(3)\times SO(3)$). Changing orientation interchanges the two bundles so we can work with either and we fix on $\Lambda^{-}$.

Write $d\theta_{i}= \psi_{i}$ and consider the linear equations for $1$-forms
$T_{i}$
\begin{equation}    \psi_{i}= T_{j}\wedge \theta_{k}- T_{k}\wedge \theta_{j}.\end{equation}
(Here, and below, we use the convention that $(ijk)$ runs over the cyclic permutations of $(123)$.) It is a fact that this system of linear equations has a unique solution. This fact is essentially the same as the usual characterisation of the Levi-Civita connection in that the covariant derivative on $\Lambda^{-}$ is
$$   \nabla \theta_{i}= T_{j}\otimes \theta_{k}-T_{k}\otimes \theta_{j}. $$
The solution of the equations (22) is
$$  -2 T_{i}= *\psi_{i} - \theta_{j}\wedge (*\psi_{k}) + \theta_{k}\wedge (*\psi_{j}), $$
where $*$ is the $*$-operator of the metric. 
The $T_{i}$ are connection forms for $\Lambda^{-}$ in the local orthonormal trivialisation $\theta_{i}$. The components of the curvature tensor of $\Lambda^{-}$ are the forms
$$  F_{i}= dT_{i} + T_{j}\wedge T_{k}. $$
This gives a way to compute the Riemann curvature tensor which is useful in some situations, such as that below. In particular we can take the anti-self-dual components $F^{-}_{i}$ of the $F_{i}$ and express them in terms of the given basis so
$$   F^{-}_{i}= \sum_{j} W_{ij} \theta_{j},$$
 Then the matrix $W_{ij}$ represents  the anti-self-dual Weyl tensor of the Riemannian metric. It is a general fact that this is symmetric and trace-free.

A particular case of this is when the forms $\theta_{i}$ are all closed. Then the $T_{i}$ vanish and we see that $\Lambda^{-}$ is flat. This means that locally we have a hyperk\"ahler metric, although to fit with standard conventions we should change orientation, so we are considering closed forms
$\omega_{i}$. In this situation the only non-vanishing component of the Riemann curvature tensor is the anti-self-dual Weyl tensor, so we can use a basis $\theta_{i}$ to compute the curvature tensor, as above. 
 
\subsection{The Gibbons-Hawking construction}
We review this well-known construction.
We start with a positive harmonic function $f$ on a domain $\Omega$ in $\bR^{3}$ and an $S^{1}$ bundle $P$ over $\Omega$ with a connection whose curvature is $-*df$. Let $\alpha$ be the connection 1-form over $P$ and $dx_{i}$ the pull-back of the standard $1$-forms on $\bR^{3}$. Then we have
$$  d\alpha= -\sum f_{i} dx_{j} dx_{k}. $$
(We will write $f_{i}, f_{ij}$ etc. for the partial derivatives of $f$.)
Set  $$\omega_{i}= \alpha \wedge dx_{i}+ f dx_{j}\wedge dx_{k}. $$ Then 
$d\omega_{i}= -f_{i} dx_{i} dx_{j}dx_{k}+ f_{i} dx_{i}\wedge dx_{j}\wedge dx_{k}=0$ and it is clear that the forms satisfy $\omega_{i}\wedge \omega_{j}= \delta_{ij} V$ , with $V=f \alpha \wedge dx_{1}\wedge dx_{2}\wedge dx_{3}$, so we have a hyperk\"ahler structure.

\

One basic example is when $\Omega=\bR^{3}-\setminus \{0\}$ and $f=4\pi^{-1} \vert x\vert^{-1}$.
Then the manifold we construct is $\bR^{4}\setminus\{0\}$ with the flat metric.
If we identify $\bR^{4}$ with $\bC^{2}$ in the usual way, the circle action can be taken to be $(z,w)\mapsto (\lambda z, \lambda^{-1} w)$ and the map from $\bC^{2}$ minus the origin to $\bR^{3}$ given by the identification with $P$ is
$$   (z,w)\mapsto ({\rm Re} (zw), {\rm Im} (zw), \vert z\vert^{2}-\vert w\vert^{2}). $$ 
Like the metric, this map  extends smoothly over the origin but we get a fixed point of the action, corresponding to the pole of $f$.   In general if we start with a hyperk\"ahler $4$-manifold $(M,\omega_{1}, \omega_{2}, \omega_{3})$ with a circle action which is Hamiltonian with respect to the three symplectic forms then the Hamiltonians $x_{i}:M\rightarrow \bR$ define a map $\underline{x}:M\rightarrow \bR^{3}$ and we recover the structure on $M$ (at least locally) from a harmonic function with poles.

Now we want to compute the curvature tensor of such a hyperk\"ahler 4-manifold.
Set $\theta_{i}= \alpha \wedge dx_{i}- f dx_{j}\wedge dx_{k}$. Then $\theta_{i}$ form an orthonormal basis for $\Lambda^{-}$ as considered in (5.1) and
$$ d\theta_{i} = -2 f_{i} dx_{1} dx_{2}dx_{3}. $$

 One finds then that the 1-forms $T_{i}$ are 
$$T_{i}= -\frac{f_{i}}{f} \alpha + \frac{1}{f^{2}}(f_{j}dx_{k}-f_{k} dx_{j}). $$
 Computing $dT_{i}+ T_{j}\wedge T_{k}$, one finds that the curvature tensor $(W_{ij})$ in this orthonormal basis is the trace-free part of the matrix
\begin{equation} \left( \frac{f_{ij}}{f^{2}}- 3 \frac{f_{i}f_{j}}{f^{3}}\right).  \end{equation}
This can also be written as $2f$ times the trace-free part of the Hessian of the function $f^{-2}$ which checks with the fact that when $f=\vert x\vert^{-1}$ the construction yields the flat metric on $\bR^{4}$. For then $f^{-2}= \vert x\vert^{2}$, the Hessian of $f^{-2}$ is twice the identity matrix and so its trace-free part is zero.
\subsection{Cone singularities}
Now return to our cone metric on $\bR^{2}\times \bR$ and let $f$ be the Green's function $f(x)=\Gamma_{p}(x)= G(x,p)$ where $p=(1,0,0)$. Locally, away from the singular set, we can identify domains in $\bR^{2}_{\beta}\times \bR$ with domains in $\bR^{3}$ and it is clear that the construction above yields a Ricci-flat metric on an $S^{1}$ bundle $P$ over the complement of the singular set and the point $p$. Another useful way to think about this is to cut the plane along the negative real axis and identify the corresponding cut 3-space with a wedge-shaped region $U$ in standard Euclidean 3-space. We perform the Gibbons-Hawking construction in the usual way over $U$, with the pole of $f$ yielding a fixed point of the action. Then we reverse the cut we made and glue appropriate points on the boundary of the 4-manifold to get our metric with cone singularity. Either way,  the upshot is that we get a 4-manifold $\overline{P}$ with an $S^{1}$ action having a single fixed point, a map $\pi:\overline{P}\rightarrow  \bR^{2}\times \bR$ and a metric $g$ on $\overline{P}$ with a cone singularity along $\pi^{-1}(S)$. 

The metric $g$ is locally hyperk\"ahler but not globally. It has a global K\"ahler structure $\omega_{1}$ corresponding to the direction of the edge of the wedge. If we choose local structures $\omega_{2}, \omega_{3}$ then parallel transport around the singular set takes the complex form $\Theta=\omega_{2}+i\omega_{3}$ to $ e^{2\pi (\beta-1) i }\Theta$. 

Now we claim that, with this global complex structure, $\overline{P}$ can be identified with $\bC^{2}$ and the singular set $\pi^{-1}(S)$ corresponds to the complex curve $C=\{ (z,w)\in \bC^{2}: zw=1\}$. For this we begin by going back to the general Gibbons-Hawking construction with a harmonic function $f$ on a domain $\Omega$ which we suppose to be the product of a domain  in the plane $x_{1}=0$ with an interval about $0$ in the $x_{1}$ co-ordinate. Trivialise the bundle $P$ by parallel transport in the $x_{1}$ direction, so the connection $1$-form is $a=a_{2} dx_{2}+a_{3}dx_{3}$. Write $\psi$ for the angular co-ordinate on the fibres of $P$. We seek a holomorphic function $h$ on $P$, for the complex structure corresponding to $x_{1}$.
In the trivialisation this amounts to solving the equations
\begin{equation}  \frac{\partial h}{\partial x_{1}}= -i f \frac{\partial h}{\partial \psi}\ \ \ , \ \ \ \frac{\partial h}{\partial x_{2}}+\frac{\partial h}{\partial x_{3}}= (a_{2}+ia_{3}) h. \end{equation}
We look for a solution 
 which has weight $1$ for the circle action, $h(\lambda z)= \lambda h(z)$. In this case $\frac{\partial h}{\partial \psi}=ih$ so  the first equation gives 
\begin{equation} h(x_{1}, x_{2}, x_{3},\psi) = e^{u} h(0,x_{2}, x_{3}) e^{i\psi} \end{equation}
where  
\begin{equation}  u(x_{1}, x_{2}, x_{3})= \int_{0}^{x_{1}} f(t,x_{2}, x_{3}) dt. \end{equation}
Conversely if we find a solution $h(0,x_{2}, x_{3})$ of the second equation in (24) over the slice $x_{1}=0$ and define $h$ using (25),(26) then the integrability condition for the complex structure implies that we obtain a solution of (24). In particular suppose that we are in the case   when $\frac{\partial f}{\partial x_{1}}$ vanishes on the plane $x_{1}=0$. This means that we can choose
$a_{2}, a_{3}$ to vanish on this plane. Thus, on this plane, the second equation in (24) is the ordinary Cauchy-Riemann equation. Given any holomorphic function $h_{0}(x_{2}+i x_{3})$ the formulae (25), (26) define a holomorphic function $h$ on $P$. The same discussion applies if we seek a function $\tilde{h}$ which transforms with weight $-1$.  We get another holomorphic function 
$$ \tilde{h}(x_{1}, x_{2}, x_{3},\psi) = e^{-u} h_{0}(0,x_{2}, x_{3}) e^{-i\psi}.$$
Thus we get a pair of holomorphic functions $(h,\tilde{h})$ on $P$ with
$$  \tilde{h} h = h_{0}(x_{2}+i x_{3}), $$
or in other words a holomorphic map from $P$ to $\bC^{2}$. This maps the lifts $\psi=0,\pi$ in $P$ of  the 2-dimensional domain in $\{x_{1}=0\}$  to the diagonal $\{z=w\}$ in $\bC^{2}$. It also maps  the subset in $P$ lying over any line $x_{2}=\xi_{2}, x_{3}=\xi_{3}$ to the plane curve $zw= h_{0}^{2}(\xi_{2}+i\xi_{3})$. 

We apply this discussion to the case when $f$ is the Green's function $\Gamma$ on $\bR^{2}\times \bR$. Of course we can only immediately fit in with the discussion above locally but we hope that the picture will be clear to the reader. By symmetry, the $\bR$ derivative of $\Gamma$ vanishes on the plane $s=0$ and we are in the position above. Moreover the symmetry taking $s$ to $-s$  lifts to a symmetry interchanging $h,\tilde{h}$. Of course one has to consider how the local construction above works around the pole, but this is just the same as in the model case of the ordinary Green's function on $\bR^{3}$. In terms of our usual co-ordinates $(r,\theta)$ on $\bR^{2}$ we define
$  h_{0}= 1- r^{c} e^{i\theta}$. This is holomorphic with respect to the given complex structure on the plane and vanishes at the pole of $\Gamma$.  The construction above produces global holomorphic functions $h,\tilde{h}$ on $\overline{P}$ with $h=\tilde{h}$ on a (real) 2-plane in $P$ which maps to the  plane $s=0$ in $\bR^{2}\times \bR$ as a double branched cover, branched over the origin. The functions satisfy $h\tilde{h}=1$ on the singular set.
So we get a holomophic map from $\overline{P}$ to $\bC^{2}$ taking the circle action on $P$ to the action $(z,w)\mapsto (\lambda z,\lambda^{-1} z)$ and mapping the singular set to the curve $zw=1$. The fact that $\Gamma(r,\theta,s)$ decays like $s^{-1}$ as $s\rightarrow \infty$, so its indefinite integral with respect to $s$ is unbounded, implies  that this map is bijective, by a straightforward argument. 

  We can also start from the opposite point of view with the complex manifold $\bC^{2}$ and the $\bC^{*}$-action $(z,w)\rightarrow (\lambda z,\lambda^{-1} w)$. We consider a locally-defined holomorphic $2$-form
$$\Theta= \beta(1-zw)^{\beta-1} dz dw. $$
This is preserved by the $\bC^{*}$-action and, locally, there is a holomorphic Hamiltonian map $H_{\bC}(z,w)= (1-zw)^{\beta}$. Although this is not well-defined globally the power $H_{\bC}^{1/\beta}$ is so,  and this gives the $\bR^{2}$ component of the map from $\bC^{2}$ to $\bR^{2}\times \bR$ which arises from the identification of $\bC^{2}$ with $\overline{P}$.

 It would take a little work to check that the metrics we have studied here really do give metrics with cone singularities of the kind we defined in Section 4---analysing the local representation in complex co-ordinates, but it seems to the author that this should not be hard.

There are several possible variants of this construction. For example,  we can use finite sums of Green's functions to get Ricci-flat
 K\"ahler metrics with cone singularities  on ALE spaces.
The example we have constructed above furnishes a plausible model for certain degenerations of metrics with cone singularities on compact manifolds. Consider a compact complex surface $X$ and a family of curves $D_{\epsilon}$ which converge as $\epsilon\rightarrow 0$ to a singular curve $D_{0}$ with one ordinary double point at $p\in X$. Suppose there are K\"ahler-Einstein metrics $\omega_{\epsilon}$ with fixed cone angle $\beta$ along $D_{\epsilon}$, for $\epsilon\neq 0$. We should expect that, after re-scaling small balls about $p$, the rescaled metrics converge to the Ricci-flat metric we have discussed above. Thus these kind of Ricci-flat, non-compact model solutions  should play the same role in the theory of metrics with cone singularities that the ordinary ALE spaces play in the standard theory. 

Another interesting application of these ideas is to supply models for the behaviour of the metrics around the  singular set. In particular we can study the growth of the curvature. Looking at  (23), we see that the curvature will be dominated by the Hessian of the harmonic function $f$ and from the discussion in Section 3 we see that this will typically be $O(r^{c-2})=O(r^{\beta^{-1}-2})$. Since $\beta^{-1}>1$ the curvature is, at least locally, in $L^{2}$ but if $\beta>1/2$ the curvature is unbounded. We expect that the like will hold for general K\"ahler-Einstein metrics with cone singularities.

\section{Conjectural picture}
In this final section we discuss what one might expect about the existence problem for metrics with cone singularities on a Fano manifold.   It is natural to think of such a metric as a solution of a distributional equation
\begin{equation} \Ric(\omega)= \beta \omega + 2\pi (1-\beta) [D]. \end{equation}
But in writing this equation we emphasise that we mean solutions of the kind we have defined precisely in Section 4.
  This equation can be compared with the  equation studied in the standard \lq\lq continuity method''
\begin{equation} \Ric(\omega)= \beta \omega + (1-\beta) \rho \end{equation}
where $\rho$ is a prescribed closed $(1,1)$ form representing $2\pi c_{1}(X)$.  There are good reasons for believing that the cone singularity problem will always have solutions for {\it small} positive $\beta$. In one direction, Tian and Yau established the existence of a  complete Ricci-flat K\"ahler metric on the non-compact manifold
$X\setminus D$ \cite{kn:TY} and one could expect that this is the limit of solutions $\omega_{\beta}$ as $\beta$ tends to $0$ (this idea, in the negative case, is mentioned by Mazzeo in \cite{kn:Mazzeo}). In another direction, it is known that, at least if $X$ has no holomorphic vector fields, solutions to (28) exist for small $\beta$ and one could perhaps view (27) as a limiting case. Szekelyhidi \cite{kn:GS}
introduced an invariant $R(X)$, defined to be the supremum of numbers $\mu$ such there is a K\"ahler metric $\Omega$ in the class $2\pi c_{1}(M)$ with $\Ric(\Omega)\geq \mu \Omega$, pointwise on the manifold. He showed that for any choice of $\rho$ this is also the supremum of the values $\beta\leq 1$ such that a solution of (28) exists. Further, it is known that $R(X)\geq \frac{n+1}{n} \alpha(X)$ where $\alpha(X)$ is Tian's invariant. The natural conjecture then is
\begin{conj}
There is a cone-singularity solution $\omega_{\beta}$ to (27) for any parameter $\beta$ in the interval $(0,R(X))$.
 If $R(X)<1$ there is no solution for parameters $\beta$ in the interval $(R(X),1)$ 
\end{conj}

Note that if $\beta=\nu^{-1}$ for an integer $\nu$, our metrics with cone singularities are {\it orbifold} metrics, so a great deal of standard theory can be brought to bear. See the recent work \cite{kn:RT} of Ross and Thomas, for example. 

\

Suppose that we are in a case when solutions   exist for small cone angles but not for cone angles close to $1$. We would like to understand how the solutions can break down at some critical cone angle. This leads into a large discussion involving notions of \lq\lq stability'' which we only want to touch on here. Recall  that in the established theory one defines the {\it Futaki invariant} of a K\"ahler manifold $Y$ with a fixed circle action. One definition is to take any invariant metric in the K\"ahler class and then set
\begin{equation}  {\rm Fut}(Y)= \int_{Y} (S- \hat{S}) H \end{equation}
where $S$ is the scalar curvature, $\hat{S}$ is the average value of the scalar curvature and $H$ is the Hamiltonian of the circle action. The key point is that in fact the Futaki invariant does not depend on the choice of metric, in a fixed K\"ahler class. There are other definitions which generalise to singular spaces and schemes. What is visible from the formula (29)  is that if $Y$ admits an invariant metric of constant scalar curvature, in the given class,   then the Futaki invariant vanishes, since in that case $S=\hat{S}$. 

Now let $\Delta\subset Y$ be a divisor invariant under the circle action and $0<\beta\leq 1$. We define a Futaki invariant of the data by
\begin{equation}  {\rm Fut}(Y,\Delta, \beta)= {\rm Fut}(Y) - (1-\beta)\left( \int_{\Delta} H - \frac{\Vol(\Delta)}{\Vol(X)} \int_{X} H \right). \end{equation}
This definition can be motivated, in the framework of metrics with cone singularity along $\Delta$,  by adding a suitable distributional term to the scalar curvature, in the manner of (27) and substituting into (29). Under plausible assumptions about the behaviour around the singular set, the definition implies that if there is an invariant constant scalar curvature metric with  cone angle $\beta$ along $D$ then ${\rm Fut}(Y,\Delta,\beta)=0$.  In particular this should apply in the K\"ahler-Einstein situation. 
\begin{conj} Let $X$ be a Fano manifold and $D$  a smooth divisor in $-K_{X}$. Suppose $\beta_{0}\leq 1$ and there are K\"ahler-Einstein metrics with cone angle $\beta$ along $D$ for $\beta<\beta_{0}$ but not for cone angle $\beta_{0}$.
Then the pair $(X,D)$ can be degenerated to a pair $(Y,\Delta)$, which has an $S^{1}$ action, and ${\rm Fut}(Y,\Delta, \beta_{0})=0$.
\end{conj}
This conjecture really needs to be fleshed out. In one direction, we should discuss  pairs $(Y,\Delta)$ with singularities. In another direction, what is really relevant is that the Futaki invariant ${\rm Fut}(Y,\Delta, \gamma)$ {\it decreases} to $0$ as $\beta$ increases to $\beta_{0}$ where the sign of $H$ is linked to the degeneration of $(X,D)$ to $(Y,\Delta)$. But the statement conveys the general idea. 

We can illustrate this, albeit still at the conjectural level, by considering  two rational surfaces $X_{1}, X_{2}$: the blow-ups of $\bC\bP^{2}$ in one or two points respectively. It is well-known that these do not admit K\"ahler-Einstein metrics and we will see that the calculation of certain Futaki invariants reproduces explicit known values of the invariants $R(X_{i})$ obtained by Szekelyhidi \cite{kn:GS} and Chi Li \cite{kn:CL}.

We begin with the case of $X_{2}$, which take to be the blow-up of $\bC\bP^{2}$ at the points $p=[1,0,0], q=[0,1,0]$. We take a smooth cubic $C$ in $\bC\bP^{2}$ through these two points, so the proper transform $D$ of $C$ is a canonical divisor  in $X_{2}$. In this case the degeneration of the pair $(X_{2}, D)$ will only involve $D$, so $Y=X_{2}$. To obtain $\Delta$ we consider the $\bC^{*}$-action on $X_{2}$ induced by $[u,v,w]\mapsto [\lambda u, \lambda v, w]$ on $\bC\bP^{2}$. We define $\Delta$ to be the limit  of $D$ under the action as $\lambda \rightarrow 0$. This is the proper transform of a singular curve
$C'$ in $\bC\bP^{2}$ which is the union of three lines through $r=[0,0,1]$ (the lines $\overline{pr}, \overline{qr}$ and one other line). We take the circle action on $Y=X_{2}$ to be the obvious one defined by the above $\bC^{*}$-action. It is then a straightforward exercise to compute the Futaki invariant
${\rm Fut}(X_{2}, \Delta, \beta)$ as a function of $\beta$. To fix signs and constants we take the Hamiltonian $H$ to vanish at $r$ and to take the value $3$ on the line at infinity $\{[u,v,0]\}$. The calculation of 
${\rm Fut}(X_{2})$ is easiest using a toric description.  One finds that
$$  {\rm Fut}(X_{2})= -2/3 . $$
Likewise
$$ \Vol(X_{2})= 7/2  \ , \Vol(\Delta)= 7$$
$$  \int_{X_{2}} H = 19/3, \int_{\Delta} H =17/2 $$
Thus  ${\rm Fut}(X_{2}, \Delta, \beta)=-\frac{2}{3}+\frac{ 25(1-\beta)}{6} $ and this vanishes when $1-\beta=4/25$.  This fits in with the result of Chi-Li that $R(X_{2})= 21/25$. In fact this is not too surprising because the calculation in \cite{kn:CL} involves essentially the same ingredients, but from a different point of view.

Notice that this discussion ties in with that in Section 5 because the curve $\Delta$ is singular. We expect that re-scaling the metrics for parameters $\beta<\beta_{0}$ around the point $r$ we will get a limit which is a  Ricci-flat metric on  $\bC^{2}$ with cone angle $\beta_{0}$ along the affine part of $C$.

There is a similar discussion for $X_{1}$. We define this to be the blow-up of $\bC\bP^{2}$ at $r$ and now we take $C$ to be a smooth cubic through $r$.
This time we degenerate in the opposite direction, taking $\lambda\rightarrow \infty$. The limit $\Delta$ is a divisor which is the sum of the proper transform of a line through $r$ and the line at infinity, taken with multiplicity $2$. With $H$ normalised to be equal to $1$ on the exceptional divisor and $3$ on the proper transform of the line at infinity, the calculation now yields
$$ {\rm Fut}(X_{1},\Delta, \beta)= \frac{2}{3}- \frac{14 (1-\beta)}{3}$$ 

and we find the critical value $\beta_{0}= 6/7$, agreeing with \cite{kn:GS}, \cite{kn:CL}. (The fact that the coefficient of $(1-\beta)$ has different signs in the two cases is connected with the fact that we take limits in opposite directions $\lambda\rightarrow 0,\infty$.) The expected behaviour of the K\"ahler-Einstein metrics as $\beta\rightarrow \beta_{0}$ is less clear in this case and we leave that discussion for another place. 

{\bf Remark}

Towards the end of the writing of this paper, preprints by Berman \cite{kn:Ber} and Chi Li \cite{kn:CL2} appeared. These both seem to be very relevant to the discussion of this section, and  to give additional evidence for the conjectural picture above.
 


\end{document}